\documentclass[paper]{geophysics}
\pdfoutput=1
\usepackage{amsmath}
\usepackage{amssymb}
\usepackage{amsfonts}
\usepackage{color}
\usepackage{array}
\usepackage{seg}
\usepackage{hyperref}

\begin{document}

\renewcommand{\figdir}{.} 

\title{Total-variation minimization with bound constraints}
\author{Musa Maharramov \href{mailto:maharram@stanford.edu}{\tt <maharram@stanford.edu>} \\
        {\&}\\
 Stewart A. Levin \href{mailto:stew@sep.stanford.edu}{\tt <stew@sep.stanford.edu>} }

\righthead{Total variation with bound constraints}
\lefthead{Maharramov and Levin}
\footer{SEP--158}
\maketitle

\begin{abstract}
We present a powerful and easy-to-implement algorithm for solving constrained optimization problems that involve $L_1$/total-variation regularization terms, and both equality and inequality constraints. We discuss the relationship of our method to earlier works of \cite{split} and \cite{cw}, and demonstrate that our approach is a combination of the augmented Lagrangian method with splitting and model projection. We test the method on a geomechanical problem and invert highly compartmentalized pressure change from noisy surface uplift observations. We conclude the paper with a discussion of possible extension to a wide class of regularized optimization problems with bound and equality constraints.     
\end{abstract}

\section{Introduction}
The primary focus of this work is a class of least-squares fitting problems with a total-variation (TV) regularization and bound model constraints:
\begin{equation}
\begin{aligned}
& \| | \nabla \mathbf{m} | \|_1 \; + \; \frac{\alpha}{2}\|\mathbf{F}(\mathbf{m})-\mathbf{d} \|^2_2 \;\rightarrow\; \min,\\
&\mathbf{m}_1\;\le \; \mathbf{m} \; \le \; \mathbf{m}_2.
\end{aligned}
\label{eq:opt}
\end{equation}
In (\ref{eq:opt}) we seek a model vector $\mathbf{m}$ such that forward-modeled data $\mathbf{F}(\mathbf{m})$ match observed data $\mathbf{d}$ in the least squares sense, while imposing blockiness-promoting total-variation (TV) regularization \cite[]{ROF} and lower ($\mathbf{m}_1$) and upper ($\mathbf{m}_2$) model bounds. Rather than using a regularization parameter as a coefficient of the regularization term, we use a data-fitting weight $\alpha$. TV regularization (also know as the Rudin-Osher-Fatemi, or ROF, model) acts as a form of model styling that helps to preserve sharp contrasts and boundaries in the model, even when spectral content of input data has limited resolution. Examples of successful geophysical application of unconstrained TV-regularized optimization are included in \cite[]{musasep1581,musasep1584,yinbin1581,yinbin1582}. The regularization provided by bounded total-variation sometimes produces sufficient smoothing side-effect on the inverted model that obviates explicit bound constraints. However, many applications still require the imposition of additional constraints regardless of the regularization. For example, reservoir pore-pressure inversion problems often come with \emph{a priori} bounds on the estimated pore pressure change, such as the pore pressure change being non-negative as a result of fluid injection (lower bound) or never exceeding a hydraulic fracturing pressure (upper bound). An example of such inversion for an unconventional reservoir from field tilt measurements is provided by \cite{musasep1583}. TV regularization is a key tool in imaging and de-noising applications \cite[]{ROF,Chambolle1997,split,cw} and require an efficient mechanism for including \emph{a priori} model constraints that can significantly reduce model space \cite[]{cw}.
While barrier or penalty function methods, such as nonlinear interior-point methods \cite[]{Nocedal}, can be used to tackle the general constrained formulation (\ref{eq:opt}), the presence of a non-differentiable $L_1$-norm total-variation term and non-quadratic penalty terms pose considerable challenges to practical implementation. A log-barrier function such as 
\begin{equation}
-\, \mathrm{const}\, \times  \sum_{i=1}^n {\log \frac{ m^i_2 - m^i } {\delta } + \log \frac {m^i - m_1^i}{\delta}},
\label{eq:penalty}
\end{equation}
where $n$ is the model space dimension, can be added to the right-hand side of the objective function to keep solution iterates away from the rectangular bounds. However, this adds a non-quadratic term to the objective function. For large-scale inversion problems with $n > 10^5$ (such as typical in geophysical applications) often only iterative gradient-based solution techniques like the nonlinear conjugate gradients \cite[]{Nocedal} are available, and adding non-quadratic terms may significantly affect convergence properties. Note that this is in addition to the challenges associated with handling the non-differentiable TV-regularization term.

\cite{cw} used a splitting approach to decouple the TV-regularized problem from enforcing the constraints. In our approach, we perform three-way splitting of problem (\ref{eq:opt}) into a smooth optimization, gradient thresholding and projection steps using the Alternating Direction Method of Multipliers (ADMM) \cite[]{ADMM}. For unconstrained TV-regularized problems this approach is equivalent to the split-Bregman method of \cite{split}. However, we integrate the projection step associated with enforcing the bound constraints into the TV-minimization loop and avoid unnecessary iterations in the minimization of a proximal term \cite[]{proximal} associated with the projection.

\section{Method}

First, we recast the TV-regularization part of (\ref{eq:opt}) as a constrained optimization problem following the approach of \cite{split}. We introduce an auxiliary variable $\mathbf{x}$ and operator $\boldsymbol\Phi: \mathbf{m}\to\mathbf{x}$ such that for isotropic TV regularization we have a vector of the model-space dimension
\begin{equation}
\boldsymbol \Phi(\mathbf{m}) \;=\; \sqrt { \left(\nabla_x\mathbf{m}\right)^2 + \left(\nabla_y{\mathbf{m}}\right)^2 },
\label{eq:PHI1}
\end{equation}
and for anisotropic regularization a vector twice the model-space dimension 
\begin{equation}
\boldsymbol \Phi(\mathbf{m}) \;=\; \begin{bmatrix} \nabla_x \mathbf{m} \\ \nabla_y \mathbf{m} \end{bmatrix}.
\label{eq:PHI2}
\end{equation}
Problem (\ref{eq:opt}) can now be reformulated with an additional equality constraint:
\begin{equation}
\begin{aligned}
& \| \mathbf{x} \|_1 \; + \; \frac{\alpha}{2}\|\mathbf{F}(\mathbf{m})-\mathbf{d} \|^2_2 \;  \;\rightarrow\; \min,\\
&\mathbf{x}\;=\; \boldsymbol \Phi(\mathbf{m}),\\
& \mathbf{m}_1\;\le \; \mathbf{m} \; \le \; \mathbf{m}_2.
\end{aligned}
\label{eq:opt2}
\end{equation}
Problem (\ref{eq:opt2}) is still a bound-constrained problem. Introducing the projection operator
\begin{equation}
\boldsymbol \Pi (\mathbf{m})\;=\; \max\{\min\{ \mathbf{m}, \mathbf{m}_2\}, \mathbf{m}_1\},
\label{eq:proj}
\end{equation}
where $\min$ and $\max$ are applied component-wise, we reduce (\ref{eq:opt2}) to a fully equality-constrained formulation:
\begin{equation}
\begin{aligned}
& \| \mathbf{x} \|_1 \; + \; \frac{\alpha}{2}\|\mathbf{F}(\mathbf{m})-\mathbf{d} \|^2_2 \;  \;\rightarrow\; \min,\\
&\mathbf{x}\;=\; \boldsymbol \Phi(\mathbf{m}),\\
& \mathbf{m}\;=\; \mathbf{y},\\
&\mathbf{y} = \boldsymbol \Pi(\mathbf{m}).
\end{aligned}
\label{eq:opt22}
\end{equation}
Following the augmented Lagrangian recipe for (\ref{eq:opt22}) while assuming the last equality constraint still enforced explicitly, we obtain a sequence of problems \cite[]{Nocedal}
\begin{equation}
\begin{aligned}
& ( \mathbf{x}^{k+1}, \mathbf{m}^{k+1})\;=\; \text{argmin }   \| \mathbf{x} \|_1 \; + \; \frac{\alpha}{2}\|\mathbf{F}(\mathbf{m})-\mathbf{d} \|^2_2 \; + \\
&   \frac{\lambda}{2} \| \mathbf{x}\;-\; \boldsymbol \Phi(\mathbf{m}) \|_2^2 - {\boldsymbol\mu_k}^T\left( \mathbf{x}\;-\; \boldsymbol \Phi(\mathbf{m}) \right)  +
\frac{\delta}{2}\|\mathbf{m} - \mathbf{y}\|^2_2 - {\boldsymbol \nu_k}^T \left(\mathbf{m}-\mathbf{y}\right) \;\rightarrow\; \min,\\
& \boldsymbol \mu_{k+1}\;=\; \boldsymbol \mu_k-\lambda \left[\mathbf{x}^{k+1}-\boldsymbol \Phi(\mathbf{m}^{k+1})\right],\\ 
& \boldsymbol \nu_{k+1}\;=\; \boldsymbol \nu_k-\delta \left[\mathbf{m}^{k+1}-\mathbf{y}\right],\; k=0,1,2,\ldots\\ 
\end{aligned}
\label{eq:opt3}
\end{equation}
Coefficients $\lambda$ and $\delta$ are any positive constants above certain problem-specific ``threshold'' values \cite[]{Nocedal}, and can be selected experimentally. Vectors $\boldsymbol \mu_k$ and $\boldsymbol \nu_k$ are vectors of multipliers that converge to the set of Lagrange multipliers for the first two equality constraints of problem (\ref{eq:opt22}). At each step, (\ref{eq:opt3}) solves an $L_1$-regularized problem with respect to the combined model vector $(\mathbf{x},\mathbf{m})$. Introducing new scaled multiplier vectors
\begin{equation}
\mathbf{b}^k\;=\; \frac{ \boldsymbol{\mu}_k}{\lambda},\; \mathbf{c}^k\;=\;\frac{\boldsymbol{\nu}_k}{\delta},\; k=0,1,2,\ldots
\label{eq:b}
\end{equation}
a little algebra shows that (\ref{eq:opt3}) is equivalent to
\begin{equation}
\begin{aligned}
 ( \mathbf{x}^{k+1},&\mathbf{m}^{k+1})\;=\; \text{argmin }   \| \mathbf{x} \|_1 \; + \; \frac{\alpha}{2}\|\mathbf{F}(\mathbf{m})-\mathbf{d} \|^2_2 \; + \\
 &   \frac{\lambda}{2} \| \mathbf{x}\;-\; \boldsymbol \Phi(\mathbf{m})- \mathbf{b}^k \|_2^2 +
 \frac{\delta}{2} \| \mathbf{m}\;-\; \mathbf{y}- \mathbf{c}^k \|_2^2 \;\rightarrow\; \min,\\
& \mathbf{b}^{k+1}\;=\; \mathbf{b}^k +  \Phi(\mathbf{m}^{k+1}) - \mathbf{x}^{k+1},\\ 
& \mathbf{c}^{k+1}\;=\; \mathbf{c}^k +  \mathbf{y}-\mathbf{m}^{k+1},\; k=0,1,2,\ldots 
\end{aligned}
\label{eq:opt4}
\end{equation}
Here we used the fact that adding a constant term $\lambda/2\|\mathbf{b}^k\|_2^2+\delta/2\|\mathbf{c}^k\|^2_2$ to the objective function obviously does not change the minimizing solution.

Problem (\ref{eq:opt22}) can be solved by iteratively projecting the current model vector $\mathbf{m}$ onto $\mathbf{y}$, then conducting the iterations (\ref{eq:opt4}) to convergence, then repeating the process. However, presence of the proximal term $\delta/2\|\mathbf{m}-\mathbf{y}\|^2_2$ in (\ref{eq:opt3}) due to the constraint $\mathbf{m}=\mathbf{y}$ means that a very accurate solution of (\ref{eq:opt4}) at early iterations is wasteful and unnecessary. We instead carry out a single iteration of (\ref{eq:opt4}) followed by the model projection:
\begin{equation}
\begin{aligned}
 ( \mathbf{x}^{k+1},&\mathbf{m}^{k+1})\;=\; \text{argmin }   \| \mathbf{x} \|_1 \; + \; \frac{\alpha}{2}\|\mathbf{F}(\mathbf{m})-\mathbf{d} \|^2_2 \; + \\
 &   \frac{\lambda}{2} \| \mathbf{x}\;-\; \boldsymbol \Phi(\mathbf{m})- \mathbf{b}^k \|_2^2 +
 \frac{\delta}{2} \| \mathbf{m}\;-\; \mathbf{y}^{k}- \mathbf{c}^k \|_2^2 \;\rightarrow\; \min,\\
& \mathbf{b}^{k+1}\;=\; \mathbf{b}^k +  \Phi(\mathbf{m}^{k+1}) - \mathbf{x}^{k+1},\\ 
& \mathbf{c}^{k+1}\;=\; \mathbf{c}^k +  \mathbf{y}^{k}-\mathbf{m}^{k+1},\\
& \mathbf{y}^{k+1}\;=\; \boldsymbol \Pi (\mathbf{m}^{k+1})\;=\; \max\{\min\{ \mathbf{m}^{k+1}, \mathbf{m}_2\}, \mathbf{m}_1\},\; k=0,1,2,\ldots\\ 
\end{aligned}
\label{eq:opt5}
\end{equation}
The iterative process (\ref{eq:opt5}) still requires soling an $L_1$-regularized problem. However, the $L_1$-norm term now involves only the vector $\mathbf{x}$. Therefore, we apply splitting, minimizing 
\begin{equation}
\begin{aligned}
\| \mathbf{x} \|_1 \; + \; \frac{\alpha}{2}\|\mathbf{F}(\mathbf{m})-\mathbf{d} \|^2_2 \; + \frac{\lambda}{2} \| \mathbf{x}\;-\; \boldsymbol \Phi(\mathbf{m})- \mathbf{b}^k \|_2^2 +
 \frac{\delta}{2} \| \mathbf{m}\;-\; \mathbf{y}^{k}- \mathbf{c}^k \|_2^2 
\end{aligned}
\label{eq:obj}
\end{equation}
alternately with respect to $\mathbf{m}$ and $\mathbf{x}$ in an inner loop of $N_1\ge 1$ cycles. Because the proximal constraint $\mathbf{m}=\mathbf{y}$ renders good fitting accuracy at early stages unnecessary, $N_1$ can be small. Further we note that the minimization of (\ref{eq:obj}) with respect to $\mathbf{x}$ (in a splitting step with $\mathbf{m}$ fixed) is given trivially by the ``shrinkage'' operator \cite[]{split}:
\begin{equation}
\mathbf{x}^{k+1}\;=\;\mathrm{shrink}\left\{\boldsymbol \Phi(\mathbf{m})+\mathbf{b}^k,\frac{1}{\lambda}\right\},
\label{eq:shrink1}
\end{equation}
where
\begin{equation}
\mathrm{shrink}\left\{ \mathbf{x}, \gamma \right\}\;=\; \frac{\mathbf{x}}{|\mathbf{x}|} \max\left(|\mathbf{x}|-\gamma,0\right),
\label{eq:shrink}
\end{equation}
and is effectively thresholding the model gradient. Our algorithm can be described by the following 5 steps:
\begin{enumerate}
\item[1] Initialization
\begin{equation}
\begin{aligned}
& \mathbf{m}^{0} \;=\; \text{starting guess},\\
& \mathbf{x}^0\;=\; \boldsymbol 0,\\
& \mathbf{y}^0\;=\;  \max\{\min\{ \mathbf{m}^{0}, \mathbf{m}_2\}, \mathbf{m}_1\},\\
& \mathbf{b}^0\;=\; \boldsymbol 0,\\
& \boldsymbol c^0\;=\; \boldsymbol 0,\\
\end{aligned}
\label{eq:init}
\end{equation}
\item[2] {\bf Outer loop}. Repeat steps 3-5 for $k=0,1,2,\ldots$
\item[3] {\bf Inner loop}. Iterate (\ref{eq:inner}) $N_1\ge 1$ times.
\begin{equation}
\begin{aligned}
 \mathbf{m}^{k+1}\;=\; & \text{argmin }  \frac{\lambda}{2} \| \mathbf{x}^k - \boldsymbol \Phi(\mathbf{m}) - \mathbf{b}^k \|_2^2 + \frac{\alpha}{2}\| \mathbf{F}(\mathbf{m}) - \mathbf{d} \|_2^2 + \\
 & \frac{\delta}{2} \|\mathbf{m} - \mathbf{y}^k - \boldsymbol c^{k} \|_2^2,\\
 \mathbf{x}^{k+1}\;=\; & \mathrm{shrink}\left\{\boldsymbol \Phi(\mathbf{m}^{k+1})+\mathbf{b}^k,\frac{1}{\lambda}\right\},\; \mathbf{x}^k\;=\;\mathbf{x}^{k+1},
\end{aligned}
\label{eq:inner}
\end{equation}
\item[4] Update the multipliers and project the model onto the bounding rectangle:
\begin{equation}
\begin{aligned}
& \mathbf{b}^{k+1}\;=\; \mathbf{b}^k + \boldsymbol \Phi(\mathbf{m}^{k+1})-\mathbf{x}^{k+1},\\
& \boldsymbol c^{k+1}\; =\; \boldsymbol c^{k} + \mathbf{y}^{k} - \mathbf{m}^{k+1},\\
& \mathbf{y}^{k+1}\;=\; \max\{\min\{ \mathbf{m}^{k+1}, \mathbf{m}_2\}, \mathbf{m}_1\}.
\end{aligned}
\label{eq:outer}
\end{equation}
\item[5] Terminate if the target accuracy is reached
\begin{equation}
\frac{\|\mathbf{m}^{k+1}-\mathbf{m}^k\|_2}{\| \mathbf{m}^k \|} \; \le \; \text{target accuracy}.
\label{eq:target}
\end{equation}
or go back to step 2 otherwise.
\end{enumerate}
Optimizing (\ref{eq:inner}) with respect to $\mathbf{m}$ is in itself a large-scale optimization problem, nonlinear for a nonlinear modeling operator $\mathbf{F}$. Solving the optimization problem (\ref{eq:inner}) exactly is unnecessary because for small $k$ (i.e., at early stages of the inversion) vector $\mathbf{y}^k$ is not the true model, vector $\mathbf{x}^k$ is far from the true model gradient, and the multipliers $\mathbf{b}^k, \mathbf{x}^k$ could be far from scaled Lagrange multipliers, s. Therefore, for large-scale problems only a few steps of an iterative method like conjugate gradients need be carried out. As the solution converges to the true solution and critical sharp contrasts in the model are identified, an iterative solver can begin to take advantage of the objective function curvature information collected at previous iterations of the \emph{outer} loop, potentially leading to a significantly faster convergence. Optimal strategies for spanning iterations of nonlinear conjugate gradients across iterations of the \emph{outer} loop of our algorithm are the subject of an upcoming report.

\section{Results}

We demonstrate our method with a test problem that simulates vertical surface uplift in response to distributed dilatational sources, mathematically equivalent to surface deformation due to pore pressure change \cite[]{SEGDEF}. Our modeling operator is defined as
\begin{equation}
\mathbf{F}(\mathbf{m})\;=\; u(x),\; u(x)\;=\; \int_0^A \frac{D^3 m(\xi) d\xi  }  {\left(D^2 + (x-\xi)^2\right)^{3/2}},
\label{eq:F}
\end{equation}
where we assume that $\mathbf{m}=m(\xi), \xi\in [0,A]$ is a relative pore pressure change along a linear segment $[0,A]$ of a reservoir at a constant depth $D$, and $\mathbf{u}=u(x),x\in [0,A]$, within a proportionality factor determined by poroelastic medium properties \cite[]{musasep1583}, is the induced vertical uplift on the surface. For demonstration purposes we consider a one-dimensional model but the results trivially extend to realistic reservoir and surface geometries. Operator (\ref{eq:F}) is a smoothing operator, and recovering sharp pressure contrasts e.g. due to permeability barriers requires model ``styling'' or regularization such as blockiness-promoting ROF model. As a true model we used a highly compartmentalized pressure model of Figure~\ref{fig:true}. In our experiments, we set $D=100$m $A=2$km, and discretized both the model and data space using a 200-point uniform grid. Random Gaussian noise with $\sigma=15\%$ of the maximum clean data amplitude was added to the clean forward-modeled data to produce the noisy observations shown in Figure~\ref{fig:tvdata}.
\multiplot{2}{tvdata,true}{width=.47\columnwidth}
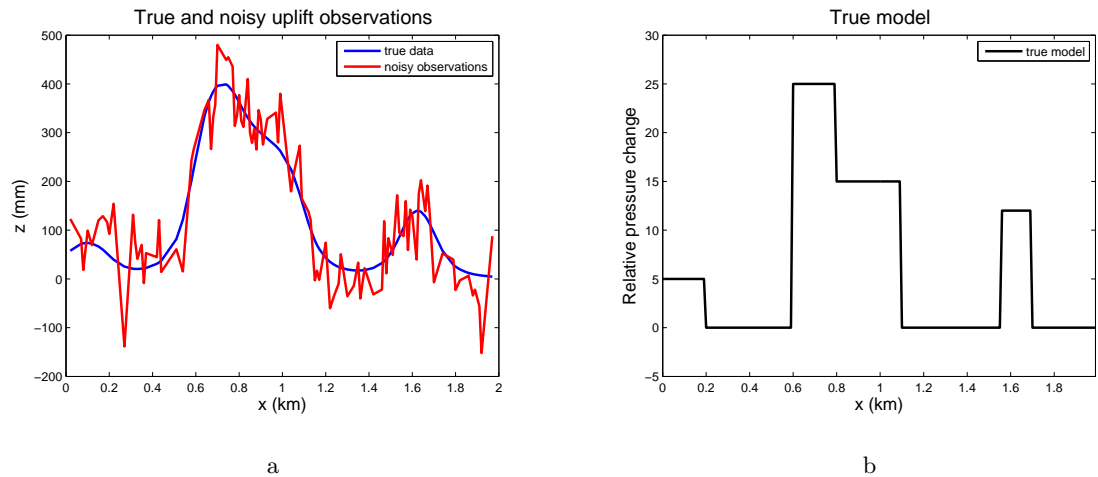
{(a) True and noisy uplift observations. Random Gaussian noise with $\sigma=15\%$ of maximum clean data amplitude was added to the clean data. (b) True model exhibibits a highly compartmentalized ``blocky'' behavior.}

The result of a TV-regularized \emph{unconstrained} inversion is shown in Figure~\ref{fig:tvinv} against the true model and a Tikhonov-regularized inversion. This result was obtained using the above algorithm by setting $\delta=0$ (no bound constraints) and using the values of $\alpha=1$ and $\lambda=2$. The TV-regularized result captures the compartmentalized picture of pressure distribution better than the highly smoothed Tikhonov regularization result. However, due to absence of bound constraints, lower pressure bounds are not honored, resulting in negative pressure areas that are not present in the true model. The result of running our bound-constrained TV-regularization algorithm is shown in Figure~\ref{fig:boundtvinv}. The imposition of bound constraints not only removed the negative relative pressure areas, but also removed the pressure spike at about $x\approx 1$km in the unconstrained inversion of Figure~\ref{fig:tvinv} that apparently had resulted from compensating negative pressures. In both the constrained and the unconstrained runs we conducted 1000 outer loop iterations with 2 inner loops cycles. However, the algorithm converged quickly, with only a few initial iterates outside a tight neighborhood of the final curve, as shown in Figure~\ref{fig:boundtvconv}.
Finally, we note that many practical implementations of bound constraints often resort to a simplistic way of enforcing the constraints: the inverted model is projected onto the bounding rectangle either once after applying a direct unconstrained solver, or at each iteration of an unconstrained solver. In this case variable $\mathbf{y}$ and the associated quadratic regularization term are not introduced into the objective function. This may result in a violation of the KKT optimality conditions where the bound constraints are active \cite[]{Nocedal}, and is demonstrated by the blue plot in Figure~\ref{fig:simplefailboundtvinv}. While the bound constraints are honored, the solution is both qualitatively and quantitatively far from optimal.

\section{Conclusions and Perspectives}

Our algorithm combines the advantages of the blockiness-promoting and edge-preserving ROF model with the ability to impose bound constraints. The splitting mechanism used for enforcing the bound constraints is naturally integrated into the split-Bregman iterations and results in no extra computational cost. The method was able to resolve compartmentalized subsurface pressure changes from noisy surface uplift observations despite the highly diffusive nature of the underlying deformation process.
\multiplot{2}{tvinv,boundtvinv}{width=.47\columnwidth}
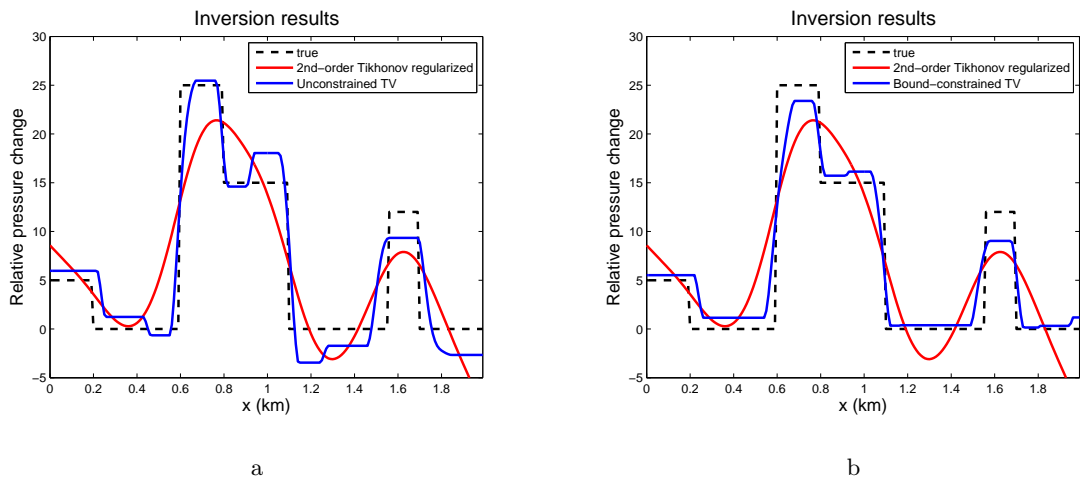
{(a) Unconstrained TV-regularized inversion. The algorithm tries to fit the data by allowing negative relative pressure changes. (b) Bound constrained TV-regularized result. Note that enforcing lower bounds resulted in a more accurate shape matching of the true model.}
The method can be implemented around any large-scale nonlinear solver such as conjugate gradients or quasi-Newton methods. Additional equality and inequality constraints can be incorporated into the algorithm using the general ADMM framework.
\multiplot{2}{simplefailboundtvinv,boundtvconv}{width=.47\columnwidth}
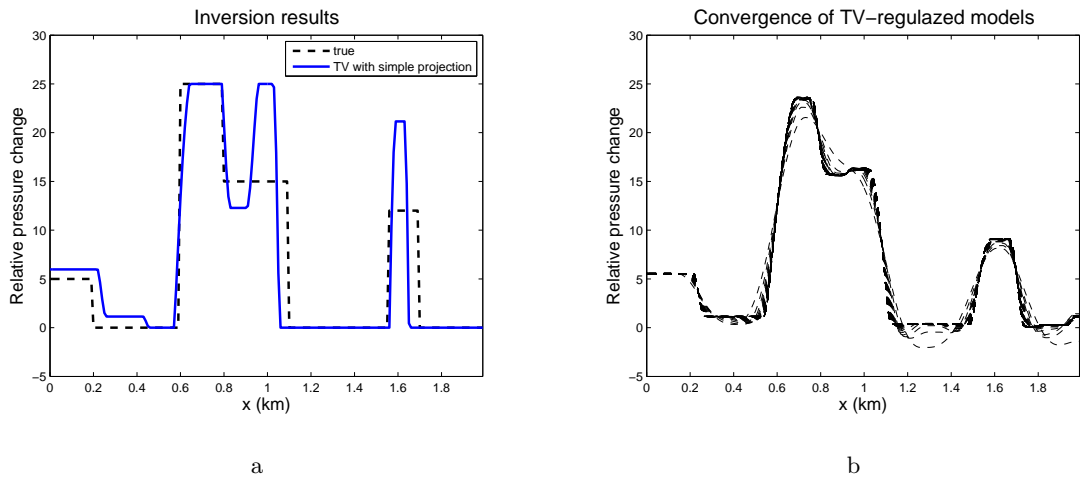
{(a) Direct imposition of the bound constraints at each iteration of the unconstrained solver resulted in a qualitatively and quantitatively wrong inversion. (b) Convergence of TV-regularized inverted models with bound constraints. The method quickly resolves both sharp contrasts and active bounds as only a few initial curves out of 1000 iterates lie outside a small neighborhood of the final curve.}

\bibliographystyle{seg}  
\bibliography{mmslbtv}

\end{document}